\documentclass[10pt,a4paper,reqno]{amsart}
\newtheorem{theorem}{\bf Theorem}[section]
\newtheorem{remark}[theorem]{Remark}

\usepackage{xcolor}
\usepackage{amsmath}
\usepackage[colorlinks=true,linkcolor=blue,citecolor=orange]{hyperref}
\newcommand{\rvs}[1]{#1}
\newcommand{\dF}[1]{\hat{#1}}
\newcommand{\newgamma}{\kappa}

\begin{document}

\title[Inverse source problem for discrete Helmholtz equation]{Inverse source problem for discrete Helmholtz equation}

\author{Roman Novikov}
\address{Roman G. Novikov, CMAP, CNRS, \'{E}cole polytechnique, Institut Polytechnique de Paris, 91128 Palaiseau, France \newline
\& IEPT RAS, 117997 Moscow, Russia}
\email{novikov@cmap.polytechnique.fr}

\author{Basant Lal Sharma}
\address{Basant Lal Sharma, Department of Mechanical Engineering, Indian Institute of Technology Kanpur, Kanpur, 208016 UP, India}
\email[Corresponding author]{bls@iitk.ac.in}

\date{\today}

\begin{abstract}
We consider multi-frequency inverse source problem for the discrete Helmholtz operator on the square lattice $\mathbb{Z}^d$, $d \ge 1$.
We consider this problem for the cases with and without phase information.
We prove uniqueness results and present examples of non-uniqueness for this problem for the case of compactly supported source function,
\rvs{and a Lipshitz stability estimate 
for the phased case is established.}
Relations with inverse scattering problem for the discrete Schr\"{o}dinger operators in the Born approximation are also 
provided.

\bigskip

\end{abstract}

\maketitle 

\section{Introduction}

We consider the discrete Helmholtz equation
\begin{equation}
\Delta \psi(x)-\lambda\psi(x)=f(x),\quad x\in\mathbb{Z}^d, d\ge 1,
\label{discschro}
\end{equation}
where
$\Delta$ is the discrete Laplacian defined by
\begin{equation}
\Delta\psi\left(x\right) =\sum_{|x^{\prime}-x|=1}
\psi\left(x^{\prime}\right), \quad x, x'\in\mathbb{Z}^d,
\label{discLap}
\end{equation}
and $f$ is a scalar source term
with sufficient decay at infinity. 
For example, our considerations 
are focused on the case when 
\begin{equation}
\text{supp } f \text{ is compact}.
\label{sufffcomp}
\end{equation}
We assume that
\begin{equation}
\lambda\in S:=
\left[ -2d,2d\right] \backslash S_{0},
\label{defS}
\end{equation}
where
\begin{eqnarray}
S_{0}&:=&\{ \pm4n\textrm{ when }d\textrm{ is even,}\nonumber\\
&&\qquad\pm2(2n+1)\textrm{ when }d\textrm{ is odd, }n\in\mathbb{Z}\textrm{ and }2n\leq d\}.
\label{defS0}
\end{eqnarray}

Discrete Helmholtz equations are closely related to discrete Schr\"{o}dinger equations which appear naturally in the tight-binding model of the electrons in crystals \cite{Bloch,Slater,Harrison}. Similar equations also appear in case of studies involving time harmonic elastic waves in lattice models of crystals \cite{Brillouin,Maradudin,Lifshitz}; see, for example, \cite{Bls0,Bls1} specially in the case $d=2$. 
In the present work, 
we consider the simplest cubic lattice and the simplest Laplacian $\Delta$ defined by \eqref{discLap}.
But our considerations can be extended to more complicated lattices and Laplacians.

For Eq. \eqref{discschro}, we consider the following solutions:
\begin{equation}
\psi^{\pm}\left( x\right)=\lim_{\nu\rightarrow\lambda\pm i0}\left(\mathcal{R}_{_{\nu}}f\right)(x),\quad x\in\mathbb{Z}^d,
\label{psipm}
\end{equation}
where
\begin{equation}
\mathcal{R}_{_{\nu}}=\left( \Delta-\nu\right) ^{-1}:l^{2}\left(
\mathbb{Z}^d\right) \longrightarrow l^{2}\left( \mathbb{Z}^d\right)
,\quad\operatorname{Im}\nu\neq0.
\label{asym1}
\end{equation}

To recall some properties of $\psi^{\pm}$, we consider the surface
\begin{equation}
\Gamma(\lambda)=\{k:k\in T^{d},\textrm{ }\phi\left(k\right) =\lambda\},\quad \lambda\in\left[ -2d,2d\right],
\label{defGamma}
\end{equation}
where
\begin{equation}
\phi\left(k\right) =2{{\sum_{i=1}^{d}}}\cos k_{i}, \quad k=(k_{1}, \dotsc, k_{d}), \quad T^{d}=\mathbb{R}^{d}/2\pi\mathbb{Z}^d.
\label{defphi}
\end{equation}

In particular,
if 
\begin{equation}
\lambda\in\Lambda=\{\zeta\in\mathbb{R}: 2d-4<\left| {{\zeta}}\right| <2d\},
\label{assumeE1}
\end{equation}
then 
\rvs{\begin{equation}
\text{$\Gamma({{\lambda}})$ is smooth strictly convex with non-zero principal curvatures, }
\label{defconvx}
\end{equation}
}
and there is a unique point
\begin{equation}
{\newgamma}={\newgamma}\left( \omega, {{\lambda}}\right) \in\Gamma(\lambda),\quad \omega\in \mathbb{S}^{d-1},
\label{defkout}
\end{equation}
where 
\begin{equation}
\frac{\nabla\phi(\newgamma)}{|\nabla\phi(\newgamma)|}=\omega.
\label{newnormal}
\end{equation}
\rvs{In connection with \eqref{assumeE1}-\eqref{newnormal}, see \cite{Shaban,Isozaki1,Isozaki2a}.}
\begin{remark}
One can consider $\Gamma(\lambda)$ to be symmetric with respect to the origin $\mathtt{O}$ in $\mathbb{R}^d$ for $\lambda>0$,
and to be symmetric with respect to the point $\mathtt{O}_\pi$ in $\mathbb{R}^d$ for $\lambda<0$, where all coordinates of $\mathtt{O}_\pi$ are equal to $\pi$.
In this case, when $2d-4<\lambda<2d$ ($-2d <\lambda<-2d+4$), the surface $\Gamma(\lambda)$ is located strictly inside the cube $[-\pi, \pi]^d$ ($[0, 2\pi]^d$, respectively); see Lemma 2 in \cite{Shaban}, and also Remark 1.1 in \cite{NovBls}.
\label{remark1p1fromfirst}
\end{remark}

Under conditions \eqref{sufffcomp}, \eqref{assumeE1}, the following asymptotic formula holds (see, for example, \cite{Shaban}):
\begin{align}
\psi^{\pm}\left( x\right) &=\frac{e^{\pm i\mu\left( \omega,\lambda\right) \left| x\right| }}{\left| x\right| ^{\frac{d-1}2}}a^{\pm}\left( \omega,\lambda\right)+O\left( \frac1{\left| x\right| ^{\frac{d+1}2}}\right),
\label{asym}\\
\omega&=\hat{x}:=x/|x|,\quad\left| {{x}}\right| \longrightarrow\infty,\quad x\in\mathbb{Z}^d,\notag
\end{align}
where
\begin{equation}
\mu\left(\omega, {{\lambda}}\right)={\newgamma}(\omega, {{\lambda}})\cdot\omega,
\label{defmu}
\end{equation}
$a^{\pm}\left( \omega,\lambda\right)$ are smooth, and the remainder can be estimated uniformly in $\omega\in\mathbb{S}^{d-1}$.

\begin{remark}

The far-field amplitudes
$a^\pm\left( \omega,\lambda\right)$ arise in \eqref{asym} for $\omega\in\Omega$, where $\Omega$ is the countable, everywhere dense subset of $\mathbb{S}^{d-1}$,
defined by
\begin{equation}
\Omega=\{\theta: \theta=x/|x|\text{ for some } x\in\mathbb{Z}^d\}.
\label{defOmg}
\end{equation}
However, we can consider $a^\pm\left( \omega,\lambda\right)$ for $\omega\in\mathbb{S}^{d-1}$ taking into account that $a^\pm$ are continuous, and even analytic, at least, under assumptions \eqref{sufffcomp}, \eqref{assumeE1}.

\end{remark}

\begin{remark}\rvs{
Eq. \eqref{discschro} can be re-written as
\begin{equation}
-(\Delta-2d) \psi(x)-(2d-\lambda)\psi(x)=-f(x),\quad x\in\mathbb{Z}^d, d\ge 1,
\label{discschronew}
\end{equation}
where $\Delta$ is the discrete Laplacian defined by \eqref{discLap}.
In addition, for the continuous case, the related equation 
can be written as
\begin{align}
&-\triangle \Psi(x) - E\Psi(x) =F(x), \ \ x\in \mathbb{ R}^d, \ d\geq 1, \ E>0,
\label{Schrodcont} 
\end{align}
where $\triangle$ is the usual continuous Laplacian and $F$ is the source function.

Therefore, for $\lambda$ close to $2d$, under assumption \eqref{defS},
the solution $\psi^-$ in \eqref{psipm}, \eqref{asym}
corresponds to the physical (outgoing radiation) solution $\Psi^+$ of \eqref{Schrodcont} for $E$ close to $0$.

In this connection, the related limit from discrete to continuous case is justified in \cite{Bls0,Bls2,Bls31} for $d=2$; see, for more details, Section 5.3 of \cite{Bls0} and
Section 4.3 of \cite{Bls2}.
\label{remarkcontn}
}\end{remark}

In the present work we consider the inverse source problem for Eq. \eqref{discschro}
consisting in recovering the source $f$ from the far-field data \rvs{$a^-$ (or $a^+$)}, for simplicity, under assumption \eqref{assumeE1}.

We show that the far-field amplitude \rvs{$a^-$ (or $a^+$)}, given in a neighbourhood of any fixed pair $\omega, \lambda$, where $\omega\in\mathbb{S}^{d-1}$ and $\lambda\in\Lambda$ 
with $\Lambda$ defined in \eqref{assumeE1}, uniquely determines $f$ under condition \eqref{sufffcomp};
see Theorem \ref{thm1}.
\rvs{A related Lipshitz stability estimate is also given; see Theorem \ref{thmstab}.}
\rvs{However}, we show that the far-field amplitude \rvs{$a^-$ (or $a^+$)}, given on $\mathbb{S}^{d-1}$ for finite but arbitrarily large number of the spectral parameter $\lambda\in\Lambda$, fails to determine $f$ under condition \eqref{sufffcomp};
see Theorem~ \ref{thm2}.

We also give uniqueness results on the phaseless inverse source problem for the discrete Helmholtz equation \eqref{discschro} with background information;
see Theorems \ref{thm4}, \ref{thm5}, \ref{thm6}.

\rvs{Additionally}, we show similarities of the multi-frequency inverse source problem for discrete Helmholtz equations 
and 
the monochromatic inverse scattering problem for discrete Schr\"{o}dinger operators 
in the Born approximation.
Our Theorems \ref{propAp1}, \ref{propAp2},
\rvs{\ref{thmAstab}},
\ref{propAp4},
\ref{propAp5}, given in Appendix \ref{appsch} for the latter case, are similar to Theorems \ref{thm1}, \ref{thm2}, \rvs{\ref{thmstab}}, 
\ref{thm5}, \ref{thm6} mentioned above \rvs{and presented in detail in Section~ \ref{mainresult}}.

\begin{remark}
If $d\ge 3$ then the case
\begin{equation}
\left| \lambda\right|
<2d-4, \quad \lambda\in S,
\label{assumeE2}
\end{equation}
is also possible.
In this case, the asymptotic formula \eqref{asym} should be replaced by 
multi-term
formula (5) in \cite{Shaban} which involves several far-field amplitudes $a^\pm$.
In addition, formula \eqref{ampshab} in our proofs should be replaced by 
formula (13) in \cite{Shaban} for the aforementioned far-field amplitudes.
It is not difficult to extend the results of the present work to the case when $\lambda$ satisfies \eqref{assumeE2}. 
However, the considerations become rather cumbersome and are not discussed in the present article.
\label{rem13p}
\end{remark}

\begin{remark}
To our knowledge, formulas \eqref{asym}, \eqref{ampshab}, and more general formulas, mentioned in Remark \ref{rem13p},
are available in the literature under assumption \eqref{sufffcomp} only.
However, these formulas, apparently, remain valid for $f$ exponentially decaying at infinity
and also for $f$ decaying at infinity as $O(|x|^{-N})$ for sufficiently large $N$ depending on $d$. 
It seems that 
$N>d$ is already sufficient with $o(|x|^{-(d-1)/2})$ in place of $O(|x|^{-(d+1)/2})$ in \eqref{asym}, by analogy with the continuous case.
In this connection we also give Remarks \ref{rem2p7} and \ref{rem2p8}.
\label{rem1p4}
\end{remark}

Note that
many important results are given in the literature for the inverse source problem for Helmholtz equation in the continuous case;
see, for example, \cite{Devaney2}--\rvs{\cite{Bourgeois}} and references therein.
\rvs{In particular, our work is motivated by obtaining discrete counterparts, for the simplest 
model \eqref{discschro}, of results of these works.
}
To our knowledge, 
the inverse source problem for the discrete Helmholtz equation \eqref{discschro} has not been considered
\rvs{yet in the literature.}
For some discrete transport equation the inverse source problem was considered in \cite{Fokas}.
Potential applications of the inverse source problem 
for discrete Helmholtz equations, including the simplest case of Eq. \eqref{discschro} can be similar to the continuous case;
\rvs{for example, within the tight-binding model of the electrons \cite{Harrison} as well as lattice dynamics in the harmonic approximation \cite{Maradudin} in crystals.
\rvs{Potential applications also include the domain of inverse source problems in forced networks; see, for example, \cite{Caputo} and references therein.}
}

On the other hand, inverse scattering for discrete Schr\"{o}dinger equations was studied already in many works; see, for example, 
\cite{Eskina2}, \cite{Case}, \cite{Zakhariev}, \cite{Guseinov}, \cite{Isozaki1}, \cite{Ando1}, \cite{Isozaki2a}
and references therein.
\rvs{Note that the analysis of these works is rather complicated specially in dimension $d\ge2$.}
\rvs{
In this connection, we were motivated by studying inverse scattering for discrete Schr\"{o}dinger equation 
within much simpler framework of the Born approximation.
In this case our results and proofs are relatively simple since the
analysis is similar to our analysis of the inverse source problem for the discrete Helmholtz equation \eqref{discschro}.
In addition, it is known that results on inverse scattering in the Born approximation admit an extension to the non-linearized case in the framework of distorted Born approximation \cite{Chew,Novikov15}.
See Appendix \ref{appsch} for more details.
}

\rvs{Note also that the studies of phaseless inverse scattering for discrete Schr\"{o}dinger equations were initiated recently in \cite{NovBls}.}

As about the inverse scattering problem for the Schr\"{o}dinger equation in the continuous case, see, for example, \cite{Devaney,Chew,Novikov2t,Novikov15,Novikovrev,Isaev,Bourgeois} and references therein.

\rvs{The} studies on inverse scattering and inverse source problems 
for both continuous and discrete cases use results on properties of Green's function for related Helmholtz operators.
In connection with these properties for the discrete case, see, for example, \cite{Morita1,Horiguchi1,Morita72,Katsura,Shaban,Isozaki1,Isozaki2a,Bls2,Bls3,Bls6,Blssurf,Garnier}
and references therein.

The main results of the present article on the inverse source problem for Eq. \eqref{discschro} are given in detail in Section \ref{mainresult}
and are proved in Section \ref{proofsecn1}.
Similar results on monochromatic inverse scattering in discrete framework in the Born approximation are given in Appendix \ref{appsch}.

\section{Main Results}
\label{mainresult}

We start with the inverse source problem for Eq. \eqref{discschro} for the phased case, that is for the case when complex values of 
\rvs{$a^-$ (or $a^+$)} are given.

\begin{theorem}
Let $f$ 
satisfy \eqref{sufffcomp}.
Then the far-field amplitude 
\rvs{$a^-=a^-\left( \omega,\lambda\right)$ (or $a^+$)}
arising in \eqref{asym}
and given in an open nonempty neighbourhood $\mathcal{N}$ of any fixed pair $\omega, \lambda$, where $\omega\in\mathbb{S}^{d-1}$ and $\lambda\in\Lambda$ 
with $\Lambda$ defined in \eqref{assumeE1},
uniquely determines $f$.
\label{thm1}
\end{theorem}

Theorem \ref{thm1}
is proved in Section \ref{proofsecn1}.

Let ${\mathcal{F}}$ denote the discrete Fourier transform defined by
the formula
\begin{equation}
\hat{u}(k)={\mathcal{F}}{u}\left( k\right) =\left( 2\pi\right) ^{-d/2}\underset{x
\in\mathbb{Z}^d}{\sum}u\left( x\right) e^{-ik\cdot x},\quad k
\in T^{d},
\label{uFT}
\end{equation}
where $u$ is a test function on $\mathbb{Z}^d$ and $k\cdot x:=\overset{d}{\underset{i=1}{\sum}}k_{i}x_{i}$.

We recall that
\begin{equation}
u\left( x\right)={\mathcal{F}}^{-1}{\hat{u}}\left( x\right)=\left( 2\pi\right)^{-d/2}\int_{T^{d}}\hat{u}(k)e^{ik\cdot x} dk, \quad x\in\mathbb{Z}^d.
\label{uinvFT}
\end{equation}

In particular, we consider the discrete Fourier transform $\dF{f}$ of the source $f$ in \eqref{discschro}.

\begin{theorem}
Let $u$ be a compactly supported function on $\mathbb{Z}^d$
and $\dF{u}={\mathcal{F}}u$.
Let 
$f={\mathcal{F}}^{-1}\dF{f}$,
where
$\dF{f}(k)=\dF{u}(k)\prod_{j=1}^{J}(\phi(k)-\lambda_j)$, $\phi$ is defined by \eqref{defphi},
$\lambda_j$ satisfy \eqref{assumeE1}, $k\in {T}^d$, and $J$ is a positive integer.
Then
$f$ 
satisfies \eqref{sufffcomp}
and its far-field amplitude 
\rvs{$a^-=a^-\left( \omega,\lambda\right)$ (and $a^+$)}
arising in \eqref{asym} vanishes identically on $\mathbb{S}^{d-1}$ for each $\lambda=\lambda_j$, $j=1,\dotsc, J$.
\label{thm2}
\end{theorem}

\begin{remark}
Under the assumptions of Theorem \ref{thm2}, the order of zero of $\rvs{a^\pm}=\rvs{a^\pm}\left(\omega, {{\lambda}}\right)$ on $\mathbb{S}^{d-1}$ for $\lambda=\lambda_j$ is equal or greater than the multiplicity of $\lambda_j$ as a root of the polynomial $\prod_{j=1}^{J}(z-\lambda_j)$ in $z$, where we do not assume that the numbers $\lambda_j$ 
are distinct.
\end{remark}

\begin{remark}
 The function $f$ in Theorem \ref{thm2} can be also written as
 \begin{equation}
 f=\left(\prod\nolimits_{j=1}^J(\Delta - \lambda_j)\right)u.
 \end{equation}
\end{remark}

Theorem \ref{thm2}
is proved in Section \ref{proofsecn1}.

 One can see that Theorem \ref{thm1} is a uniqueness result on multi-frequency inverse source problem for Eq. \eqref{discschro}, whereas 
Theorem \ref{thm2} gives a large class of non-uniqueness examples for this problem.
However, there is no contradiction between these two theorems: the first theorem involves infinitely many 
frequencies 
in contrast to the second one.

\begin{remark}
Theorems \ref{thm1} and \ref{thm2} have analogues for the continuous case.
In this case, Eq. \eqref{discschro} is considered for $x\in\mathbb{R}^d$, where $\Delta$ is the standard Laplacian, $\lambda<0$, and $\phi(k)=-k^2, ~k\in\mathbb{R}^d$.
The analogue of Theorem \ref{thm1} for the continuous case is well-known.
In the continuous analogue of Theorem \ref{thm2}, the function $u$ should be also sufficiently smooth on $\mathbb{R}^d$ in order to ensure a regularity of the function $f$ on $\mathbb{R}^d;$
see \cite{Devaney2,Burov} for the continuous case at least at a single frequency.
\end{remark}

\begin{remark}
It is also interesting to compare 
Theorem \ref{thm2}
with 
results on
recovering point sources for the continuous Helmholtz equation at a single 
frequency.
The latter problem is uniquely solvable (see
\cite{Bao2} and references therein) in contrast with Theorem \ref{thm2}.
\end{remark}

\begin{remark}
Theorem \ref{thm1} remains valid for 
$f$ exponentially decaying at infinity,
i.e.
\begin{equation}
|f(x)|\le \beta e^{-\alpha |x|}, \quad x\in\mathbb{Z}^d,
\label{fexpdec}
\end{equation}
for some $\alpha, \beta>0$,
under the assumption that formulas
\eqref{asym} and \eqref{ampshab} remain valid; see Remark \ref{rem1p4}.
\label{rem2p7} 
\end{remark}

\begin{remark}
Theorem \ref{thm1} 
is not 
valid, in general,
for $f$ decaying at infinity as $O(|x|^{-\infty})$
even if formulas \eqref{asym}
and \eqref{ampshab} remain valid.
Counter-examples can be constructed using formulas \eqref{uinvFT} with $f$ in place of $u$, where $\dF{u}$ is infinitely smooth on ${T}^d$ and 
identically zero on 
$\bigcup_{(\theta,\zeta)\in\mathcal{N}}\newgamma(\theta, \zeta)$
mentioned in formula \eqref{eqnfanlgam} in
the proof of
Theorem \ref{thm1}.
\label{rem2p8} 
\end{remark}

\rvs{
In addition to the uniqueness theorem \ref{thm1}, we also have the following Lipshitz stability result.
\begin{theorem}
Suppose that $f_1$ and $f_2$ satisfy \eqref{sufffcomp},
and
$\mathcal{N}$ is the same as in Theorem \ref{thm1},
and
$\overline{\mathcal{N}}\subset \mathbb{S}^{d-1}\times\Lambda$.
Suppose in addition
that  $\text{supp }f_1, \text{supp }f_2 \subset D$, where $D$ is a bounded domain in $\mathbb{R}^d$.
Then the following estimate holds:
    \begin{equation}
        \lVert f_2-f_1 \rVert_{\ell_2(D\cap\mathbb{Z}^d)}\le C^\pm_{D,\mathcal{N}} \lVert a^\pm_2-a^\pm_1 \rVert_{L_2(\mathcal{N})},
    \label{eqthmstab}
    \end{equation}
    where
    $a^\pm_1$, $a^\pm_2$ are the far-field amplitudes for $f_1, f_2$, respectively, and $C^\pm_{D,\mathcal{N}}$ are positive constants depending on $D$ and $\mathcal{N}$ only.
    \label{thmstab}
\end{theorem}

Theorem \ref{thmstab}
is proved in Section \ref{proofsecn1}.
This proof involves a general idea going back,
at least, to 
\cite{Bourgeois}.
}

We also consider the inverse source problem for Eq. \eqref{discschro} for the phaseless case, that is for the case when only values of $|\rvs{a^\pm}|^2$ are given.
We deal with this case using the method of background information in phaseless inverse problems; see \cite{Perutz,
Aktosun,Novikov2t,Hohage22,Sun23} and references therein.

Like in \cite{Novikov2t} for the continuous case, 
we consider $f_0$ and $f$ on $\mathbb{Z}^d$, 
assuming that
\begin{align}
 \text{supp }f_0 \subset D_0,\quad
\text{supp }f \subset D,
\label{eqsupp1}
\end{align}
where
\begin{align}
D_0, D \text{ are open convex bounded domains in }\mathbb{R}^d,\notag\\
D_0 \cap \mathbb{Z}^d \ne \emptyset, \quad D\cap \mathbb{Z}^d \ne \emptyset.
\label{eqsupp}
\end{align}

Let 
\begin{align}
&\textrm{diam } D = \sup_{x, y \, \in D} |x-y|, \label{diamdis} \\
&- D:=\{-x: \, x \in D\}, \label{17.1dis} \\
&D_1 + D_2:=\{x+y: \, x \in D_1,\, y \in D_2\}, \label{17dis}
\end{align}	
where $D, D_1, D_2$ are bounded sets in $\mathbb{R}^d$.

\begin{theorem}
Let $f_0, f$ satisfy \eqref{eqsupp1}, \eqref{eqsupp}, where $f_0\not\equiv0$.
Then the following formulas hold:
	\begin{align}
	&{\mathcal{F}}f(p) := (\overline{{\mathcal{F}}f_0}(p))^{-1}{\mathcal{F}}q(p), \quad p \in {T}^d, \label{eq19}\\
	&q(x):= \chi_{D-D_0}(x)(u(x)-(2\pi)^{-d/2} \sum_{y \in \, D_0\cap\mathbb{Z}^d}f_0(x+y)\overline{f_0(y)}), \label{eq22p1}\\
	&u(x):={\mathcal{F}}^{-1}(|{\mathcal{F}}(f+f_0)|^2)(x)
\text{ if } \textrm{dist } (D, \, D_0) > \textrm{diam } \, D,
 \label{eq21}\\
	&u(x):={\mathcal{F}}^{-1}(|{\mathcal{F}}(f+f_0)|^2)(x)-{\mathcal{F}}^{-1}(|{\mathcal{F}}f|^2)(x)
\text{ if } \textrm{dist } (D, \, D_0) > 0 \label{eq38},
	\end{align}
where $x\in\mathbb{Z}^d.$
\label{thm4}
\end{theorem}

Note that, under our assumptions, formula \eqref{eq19} is defined correctly taking into account that Meas ${\mathcal{Z}}=0$ in ${T}^d$, where $\mathcal{Z}=\{p\in{T}^d: {\mathcal{F}}f_0(p)=0\}$.

Note also that Theorem \ref{thm4} is just a 
result on phase retrieval for the Fourier transform $\mathcal{F}$ defined by \eqref{uFT}.
For known results on phase retrieval for different Fourier transforms, see, for example, \cite{Klibanov,Barnett,Novikov2t,Hohage22} and references therein.
Using Theorem \ref{thm4}, we obtain the following results on the phaseless inverse source problem for Eq. \eqref{discschro}.

\begin{theorem}
 Let $f_0, f$ satisfy \eqref{eqsupp1}, \eqref{eqsupp}, where $f_0\not\equiv0$,
$\textrm{dist } (D, \, D_0) > \textrm{diam } \, D.$
Let $\rvs{a_1^-}$ be the far-field amplitude for Eq. \eqref{discschro}
with $f$ replaced by $f+f_0$.
Then the intensity $|\rvs{a_1^-}|^2=|\rvs{a_1^-}\left( \omega,\lambda\right)|^2$ given in an open nonempty neighbourhood of any fixed pair $\omega, \lambda$, where $\omega\in\mathbb{S}^{d-1}$ and $\lambda\in\Lambda$ 
with $\Lambda$ defined in \eqref{assumeE1},
uniquely determines $f$, assuming that $f_0$ is known a-priori.
\label{thm5}
\end{theorem}

\begin{theorem}
 Let $f_0, f$ satisfy \eqref{eqsupp1}, \eqref{eqsupp}, where $f_0\not\equiv0$,
$\textrm{dist } (D, \, D_0) > 0.$
Let $a$, $\rvs{a_1^-}$ be the far-field amplitude\rvs{s} for Eq. \eqref{discschro}
and for Eq. \eqref{discschro}
with $f$ replaced by $f+f_0$, respectively.
Then the intensities $|\rvs{a}^-|^2=|\rvs{a}^-\left( \omega,\lambda\right)|^2$ and $|\rvs{a_1^-}|^2=|\rvs{a_1^-}\left( \omega,\lambda\right)|^2$ given in an open nonempty neighbourhood of any fixed pair $\omega, \lambda$, where $\omega\in\mathbb{S}^{d-1}$ and $\lambda\in\Lambda$ 
with $\Lambda$ defined in \eqref{assumeE1},
uniquely determines $f$, assuming that $f_0$ is known a-priori.
\label{thm6}
\end{theorem}

\rvs{Theorems \ref{thm5} and \ref{thm6} are also valid with ${a^-}$, ${a_1^-}$ replaced by ${a^+}$, ${a_1^+}$.}

Theorems \ref{thm4}, \ref{thm5} and \ref{thm6} are proved in Section \ref{proofsecn2}.

\begin{remark}
 For the continuous case, the analogues of Theorems \ref{thm4}, \ref{thm5}, \ref{thm6} are given in \cite{Novikov2t}; see also \cite{Sun23} and references therein.
\end{remark}

\section{Proofs of Theorems \ref{thm1}, \ref{thm2} and \ref{thmstab}}
\label{proofsecn1}
We consider $\phi$ defined by \eqref{defphi},
$\Gamma(\lambda)$ defined by \eqref{defGamma},
and $\dF{f}={\mathcal{F}}f$ defined by \eqref{uFT}.

For compactly supported complex valued $f$,
under assumption \eqref{assumeE1}, due to formula (13) from \cite{Shaban}, we have that
\begin{equation}
\rvs{a^\pm}\left( \omega,\lambda\right) =\frac{\sqrt{2\pi}\dF{f}\left(
\newgamma\left(\rvs{\pm} \omega,\lambda\right) \right) e^{\rvs{\pm}i(\sigma +2)\frac{\pi}{4}}}
{\sqrt{\left| K(\omega,\lambda))\right| }\left| \nabla\phi(\newgamma\left(\omega,\lambda\right) )\right| },\quad \omega\in \mathbb{S}^{d-1},
\label{ampshab}
\end{equation}
where 
$\rvs{a^\pm}=\rvs{a^\pm}\left( \omega,\lambda\right)$ are the far-field amplitude\rvs{s} in \eqref{asym},
$\newgamma(\omega,\lambda)$ is the point in \eqref{defkout},
$K(\omega,\lambda)$ is the total curvature (i.e. the product of principal curvatures) of $\Gamma(\lambda)$ at the point $\newgamma\left( \omega,\lambda\right)$,
and $\sigma=d-1.$ 
\rvs{More precisely, formula \eqref{ampshab} is mentioned in \cite{Shaban} for $\rvs{a^+}$
for real valued $f$.
In addition, $a^-=\overline{a^+}$ if $f$ is real, as also mentioned in \cite{Shaban}. 
In formula \eqref{ampshab}, we also use the symmetries of $\Gamma$ mentioned in Remark \ref{remark1p1fromfirst}.

Formulas \eqref{ampshab} for complex valued $f$ follow from the aforementioned results of \cite{Shaban} for real $f$, linear dependence of $a^\pm$ on $f$, the definition of $\hat{f}$ via \eqref{uFT} and the following formulas:
$\newgamma\left(-\omega,\lambda\right)=-\newgamma\left(\omega,\lambda\right),\lambda\in\Lambda,\lambda>0;$
$\newgamma\left(-\omega,\lambda\right)-0_\pi=-\newgamma\left(\omega,\lambda\right)+0_\pi, \lambda\in\Lambda, \lambda<0;$
$e^{2i0_\pi x}=1, x\in\mathbb{Z}^d$.}

Note that if $f$ is compactly supported on $\mathbb{Z}^d$, then 
\begin{equation}
\text{$\dF{f}$ is real-analytic on ${T}^d$.}
\label{eqnfanly}
\end{equation}

Note also that if $(\omega, \lambda)$ is a fixed pair as in Theorem \ref{thm1}, 
and $\mathcal{N}$ is its open nonempty neigbourhood in $\mathbb{S}^{d-1}\times\Lambda$,
where $\Lambda$ is defined in \eqref{assumeE1},
then there is an open nonempty neigbourhood $\mathcal{K}^{\rvs{\pm}}$ of $\newgamma({\rvs{\pm}}\omega, \lambda)$ in ${T}^d$
such that
\begin{equation}
\mathcal{K}^{\rvs{\pm}}\subseteq\bigcup_{(\theta,\zeta)\in\mathcal{N}}\newgamma({\rvs{\pm}}\theta, \zeta).
\label{eqnfanlgam}
\end{equation}
\rvs{We obtain \eqref{eqnfanlgam} considering
$\Gamma(\lambda)$
as in Remark \ref{remark1p1fromfirst} and $T^d$ as the corresponding cube in the same remark.
We use that:\\
\begin{equation}
\text{$\kappa_\zeta:\mathbb{S}^{d-1}\to\Gamma(\zeta)$ is a homeomorphism,}
\end{equation}
where
$\kappa_\zeta=\kappa(\theta,\zeta)$
at fixed $\zeta\in\Lambda$;
\begin{equation}
\text{$\kappa:\mathbb{S}^{d-1}\times I_{\lambda,\epsilon}\to\Gamma_{\lambda,\epsilon}:=\cup_{\zeta\in I_{\lambda,\epsilon}}\Gamma(\zeta)$ is a homeomorphism,}
\label{eqkappa}
\end{equation}
where
$\kappa=\kappa(\theta,\zeta)$,
$I_{\lambda,\epsilon}:=(\lambda-\epsilon, \lambda+\epsilon)$ for $\lambda\in\Lambda$,
and some $\epsilon>0$ such that $I_{\lambda,\epsilon}\subset\Lambda$;
\begin{equation}
\text{$\Gamma_{\lambda,\epsilon}$ is an open domain in $T^d$,}
\label{eqkappa1}
\end{equation}
where $T^d$ is considered as a cube as in Remark \ref{remark1p1fromfirst};
\begin{equation}
\kappa(\pm\omega, \lambda)\in\Gamma(\lambda)\subset\Gamma_{\lambda,\epsilon}.
\label{eqkappa2}
\end{equation}
}

\rvs{Property \eqref{eqnfanlgam} follows from
\eqref{eqkappa}--\eqref{eqkappa2}.}

\rvs{In turn, property \eqref{eqkappa2} follows from our assumptions and the definition of $\Gamma_{\lambda,\epsilon}$,
 property \eqref{eqkappa1} follows from the definition of $\Gamma(\zeta)$ via \eqref{defGamma},\eqref{defphi}.
 Next, $\kappa_\zeta$ reduces to the inverse of Gauss map for $\Gamma(\zeta)$
by definition of $\kappa$. 
In addition, due to \eqref{defconvx}, $\Gamma(\zeta)$ is smooth and strictly convex,
and it is known that the Gauss map is a diffeomorphism for smooth, strictly convex surfaces (as recalled, for example, in \cite{Isozaki2a}).
The bijectivity of $\kappa$ in \eqref{eqkappa} follows from the 
bijectivity of $\kappa_\zeta$
and
the property that $\Gamma(\zeta_1)\cap\Gamma(\zeta_2)=\emptyset$ if $\zeta_1\ne\zeta_2$.
For homoeomorphism in \eqref{eqkappa}, one should also use smooth
dependence of $\Gamma(\zeta)$ on $\zeta\in\Lambda.$
In this connection, one can also use explicit parametrization of $\Gamma(\zeta)$ by points of $\mathbb{S}^{d-1}$ as in formula (3.15) in \cite{Isozaki1}.
}

Theorem \ref{thm1} is proved as follows.

Due to formulas \eqref{ampshab}, \eqref{eqnfanlgam},
the function $a^{\rvs{\pm}}$, given in a neigbourhood of any fixed pair $(\omega, \lambda)$ as in Theorem \ref{thm1},
determines
$\dF{f}$ in a neighbourhood of $\newgamma\left({\rvs{\pm}}\omega,\lambda\right)$ in $T^{d}$.
In turn,
in view of \eqref{eqnfanly},
$\dF{f}$ in a neighbourhood of $\newgamma\in T^{d}$
uniquely determines $\dF{f}$ on $T^{d}$ via analytic continuation, and then determines $f$ using \eqref{uinvFT}.

This completes the proof of Theorem \ref{thm1}.

Let 
\begin{equation}
\delta(x)=\begin{cases}
1\text{ for }x=0,\\
0\text{ for }x\in\mathbb{Z}^d\setminus\{0\}.
\end{cases}
\label{Kdeltafn}
\end{equation}
We also recall that 
\begin{align}
(2\pi)^{d/2}{\mathcal{F}}^{-1}(\dF{u}_1\dF{u}_2)(x)
=u_1\ast u_2(x)
=\sum_{y\in\mathbb{Z}^d}u_1(x-y) u_2(y),\quad x\in\mathbb{Z}^d,\nonumber\\
\text{where }
\dF{u}_j={\mathcal{F}} u_j, \quad j=1,2,
\label{defconv}
\end{align}
where $u_1, u_2$ are test functions on $\mathbb{Z}^d.$

Theorem \ref{thm2} is proved as follows.

In fact, it is sufficient to prove that $f$ is compactly supported on $\mathbb{Z}^d.$
The rest follows from the definition of $\dF{f}$ in Theorem \ref{thm2},
formula \eqref{ampshab}, definition of $\Gamma(\lambda)$ by formulas \eqref{defGamma},
\eqref{defphi},
and the property that $\newgamma({\rvs{\pm}}\omega,\lambda)\in\Gamma(\lambda)$ {\rvs{for $\omega\in\mathbb{S}^{d-1}$}}.

\rvs{The proof that $f$ is compactly supported is as follows.}

\rvs{Using the definition of $\hat{f}$ and formula \eqref{defconv},
we obtain
\begin{equation}
f(x)
=u\ast(f_{\lambda_J}\ast f_{\lambda_{J-1}}\ast\dotsc \ast f_{\lambda_1})(x), x\in\mathbb{Z}^d,
\label{flamcomp10}
\end{equation}
where
\begin{equation}
f_{\lambda}:={\mathcal{F}}^{-1}(\phi-\lambda).
\label{flamcomp10a}
\end{equation}
In addition, we have that}
\begin{equation}
{\mathcal{F}}^{-1}\phi(x)
=\left( 2\pi\right)^{d/2}\sum_{|x'|=1}\delta(x-x'),\quad x\in\mathbb{Z}^d,
\label{flamcomp11}
\end{equation}
\begin{equation}
{\mathcal{F}}^{-1}\lambda(x)
=\left( 2\pi\right)^{d/2}\lambda\delta(x),\quad x\in\mathbb{Z}^d,
\label{flamcomp12}
\end{equation}
\rvs{where ${\mathcal{F}}^{-1}\lambda$ is the inverse Fourier transform of a constant $\lambda.$ 

From \eqref{flamcomp10a}-\eqref{flamcomp12}, we see that}
\begin{equation}
f_{\lambda}(x)
=\left( 2\pi\right)^{d/2}(\sum_{|x'|=1}\delta(x-x')-\lambda\delta(x)),\quad x\in\mathbb{Z}^d,
\label{flamcomp13a}
\end{equation}
\rvs{and, in particular,}
\begin{equation}
\text{$f_{\lambda}$ is compactly supported on $\mathbb{Z}^d$}.
\label{flamcomp13}
\end{equation}

\rvs{The property that 
$f$ is compactly supported on $\mathbb{Z}^d$ follows from formula \eqref{flamcomp10}, 
property \eqref{flamcomp13} for $\lambda=\lambda_1, \dotsc, \lambda_J$,
the assumption that 
$u$ is compactly supported on $\mathbb{Z}^d$,
and the property that the convolution of compactly supported functions is compactly supported.}

Theorem \ref{thm2} is proved.

\rvs{Theorem \ref{thmstab} is proved as follows.

We consider the linear maps
\begin{equation}
\mathcal{T}^\pm:\ell_2(D\cap \mathbb{Z}^d)\to L_2(\mathcal{N}),
\quad \mathcal{T}^\pm f=a^\pm,
\end{equation}
where
$f$ is supported on $D\cap \mathbb{Z}^d$ and considered as $f|_{D\cap \mathbb{Z}^d}$, and $a^\pm$ are the far-field amplitudes for $f$ and considered as $a^\pm|_{\mathcal{N}}$.

Note that the linearity of $\mathcal{T}^\pm$ follows from formula \eqref{ampshab}.
The property that $a^\pm\in L_2(\mathcal{N})$ follows from \eqref{ampshab}, the assumption that $\overline{\mathcal{N}}\subset \mathbb{S}^{d-1}\times\Lambda$,
and the fact that the denominator in \eqref{ampshab} does not vanish for $\lambda\in\Lambda$.}

\rvs{We identify $\ell_2(D\cap\mathbb{Z}^d)$ with $\mathbb{C}^M,$
where
\begin{equation}
M:=|D\cap \mathbb{Z}^d|,
\end{equation}
and $|\mathcal{A}|$ denotes the number of points in $\mathcal{A}.$}

\rvs{Due to Theorem \ref{thm1}, the maps $\mathcal{T}^\pm$ are injective.
Therefore, 
\begin{equation}
\text{range $\mathcal{T}^\pm$ can be identified with $\mathbb{C}^M$ with the norm induced by $L_2(\mathcal{N})$,}
\end{equation}
for each $\mathcal{T}^+$ and $\mathcal{T}^-$.
Next, we identify $\mathcal{T}^\pm$ with $M\times M$ complex matrices with non-zero determinant.

The estimate \eqref{eqthmstab} follows from the identifications above.
}

\section{Proofs of Theorems \ref{thm4}, \ref{thm5} and \ref{thm6}}
\label{proofsecn2}
The
proof of Theorem \ref{thm4} is similar to the proof of Theorems 3.1 and 3.2 
in \cite{Novikov2t}
for the case of continuous Fourier transform.

To prove Theorem \ref{thm4}, we use formula \eqref{defconv} and the formulas
\begin{equation}
 \overline{{\mathcal{F}}u}={\mathcal{F}}\,\widetilde{u}, \quad \widetilde{u}(x)=\overline{{u}(-x)},
 \label{eq34}
\end{equation}
\begin{equation}
\left( 2\pi\right)^{-d/2}{\mathcal{F}}(u_1\ast u_2)
={\mathcal{F}}u_1\,{\mathcal{F}}u_2,
\end{equation}
where $u, u_1, u_2$ are test functions.

We have that
\begin{equation}\begin{split}
\label{Idis}
\Sigma &= (2\pi)^{d/2}{\mathcal{F}}^{-1}(|{\mathcal{F}}(f+f_0)|^2) = (f+f_0)*(\widetilde{f}+\widetilde{f_0}) \\
 & =\sum_{y\in\mathbb{Z}^d}(f(x-y)+f_0(x-y))\overline{(f(-y)+f_0(-y))}\\ 
 &	=\sum_{y\in\mathbb{Z}^d}f(x-y)\overline{f(-y)}+\sum_{y\in\mathbb{Z}^d}f_0(x-y)\overline{f(-y)}+ \\
 &+\sum_{y\in\mathbb{Z}^d}f(x-y)\overline{f_0(-y)}+\sum_{y\in\mathbb{Z}^d}f_0(x-y)\overline{f_0(-y)} 
 =: \Sigma_1+\Sigma_2+\Sigma_3+\Sigma_4,
\end{split}\end{equation}
where $x\in\mathbb{Z}^d.$

Let 
\begin{equation}\label{brdis}
B_r = \{x \in \mathbb{R}^d: |x|<r\}.
\end{equation}

Note that 
\begin{align}
& \Sigma_1(x) = \sum_{y \in -D\cap\mathbb{Z}^d}f(x-y)\overline{f(-y)},\\
&\textrm{supp }\, \Sigma_1 \subset B_r, \quad r={\textrm{diam }D} , \label{33dis}\\
&\Sigma_2(x) = \sum_{y\in -D\cap\mathbb{Z}^d}f_0(x-y)\overline{f(-y)}, \label{34.1dis} \\&\textrm{supp } \, \Sigma_2 \subset D_0 - D, \label{35dis} \\
&\Sigma_3(x) = \sum_{-y \in \, D_0\cap\mathbb{Z}^d}f(x-y)\overline{f_0(-y)}, \label{36dis}\\
&\textrm{supp } \, \Sigma_3 \subset D- D_0, \label{37dis}\\
&\Sigma_4(x) = \sum_{-y \in \, D_0\cap\mathbb{Z}^d}f_0(x-y)\overline{f_0(-y)}, \\
&\textrm{supp } \, \Sigma_4(x) \subset B_r, \quad r={\textrm{diam }\, D_0},\label{39dis}
\end{align}
where $x\in\mathbb{Z}^d,$
and $D_0-D, D-D_0, B_r$ are defined 
according to \eqref{17.1dis}, \eqref{17dis}, \eqref{brdis}.

Similar to 
\cite{Novikov2t}, we have that:
\begin{equation}\label{50dis}
(2\pi)^{-d/2}{\mathcal{F}}\Sigma_3 = {\mathcal{F}}f\,\, {\mathcal{F}}\,\widetilde{f_0}; 
\end{equation}
\begin{equation}\label{59dis}
\Sigma_3(x) = \chi_{D-D_0}(x) (\Sigma(x) - \Sigma_4(x)),\text{ if } \textrm{dist } (D, \, D_0) > \textrm{diam } \, D;
\end{equation}
\begin{equation}
\Sigma_3(x) = \chi_{D-D_0}(x)(\Sigma(x) - \Sigma_1(x) - \Sigma_4(x)),
\text{ if } \textrm{dist } (D, \, D_0) > 0.
\label{55dis}
\end{equation}
Theorem \ref{thm4} follows from formulas \eqref{Idis}, \eqref{50dis}, \eqref{59dis}, \eqref{55dis}, and \eqref{eq34}.

Theorems \ref{thm5} and \ref{thm6} are proved as follows.

For compactly supported $f$,
under assumption \eqref{assumeE1}, 
as a corollary of formula \eqref{ampshab}, we have that 
\begin{equation}
|a^\rvs{\pm}\left( \omega,\lambda\right)|^2 =\frac{{2\pi}|\dF{f}\left(
\newgamma\left(\rvs{\pm} \omega,\lambda\right) \right)|^2}
{{\left| K(\omega,\lambda))\right| }\left| \nabla\phi(\newgamma\left(\omega,\lambda\right) )\right| ^2 },\quad \omega\in \mathbb{S}^{d-1}.
\label{ampshabphase}
\end{equation}

For compactly supported $f$,
as a corollary of \eqref{eqnfanly}, we have that 
\begin{equation}
|\dF{f}|^2=\dF{f}~\overline{\dF{f}\,}\text{ is real-analytic on }{T}^d.
\label{analfmod}
\end{equation}

Due to formulas \eqref{eqnfanlgam}, \eqref{ampshabphase}, 
the function $|a^\rvs{\pm}|^2$, given in a neigbourhood of any fixed pair $(\omega, \lambda)$ as in Theorem \ref{thm6},
determines
$|\dF{f}|^2$ in a neigbourhood of $\newgamma(\rvs{\pm}\omega, \lambda)$ in ${T}^d$.
Similarly,
the function $|a_1^\rvs{\pm}|^2$ given in a neigbourhood of any fixed pair $(\omega, \lambda)$ as in Theorems \ref{thm5} and \ref{thm6},
determines $|\dF{f}+\dF{f}_0|^2$ in a neigbourhood of $\newgamma(\rvs{\pm}\omega, \lambda)$ in ${T}^d$.

In turn, in view of \eqref{analfmod},
$|\dF{f}|^2$ in a neighbourhood of $\newgamma\in T^{d}$
uniquely determines $|\dF{f}|^2$ on $T^{d}$ via analytic continuation.
Similarly, $|\dF{f}+\dF{f}_0|^2$ in a neighbourhood of $\newgamma\in T^{d}$
uniquely determines $|\dF{f}+\dF{f}_0|^2$ on $T^{d}$.

Finally, the use of Theorem \ref{thm4} completes the proof of Theorems \ref{thm5} and \ref{thm6}.

\appendix
\section{Inverse scattering for discrete Schr\"{o}dinger operators in the Born approximation}
\label{appsch}

We consider the discrete Schr\"{o}dinger equation
\begin{equation}
\Delta \psi(x)+v(x)\psi(x)-\lambda\psi(x)=0,\quad x\in\mathbb{Z}^d, d\ge 1,
\label{discschro2}
\end{equation}
where
$\Delta$ is the discrete Laplacian defined by \eqref{discLap},
$v$ is scalar potential 
such that
\begin{equation}
\textrm{supp }v\subset {D},
\label{defv2}
\end{equation}
${D}$ is an open bounded domain in $\mathbb{R}^d$, ${D}\cap\mathbb{Z}^d\ne\emptyset$, and $\lambda\in S$ as in \eqref{defS}.

For Eq. \eqref{discschro2} we consider the 
scattering 
solutions 
\begin{equation}
\psi^{\rvs{-}}(x, k)=\psi_{0}+\psi_{sc}^{\rvs{-}},
\label{psiplus2}
\end{equation}
where 
\begin{equation}
\psi_{0}(x, k)=e^{i {{k}} \cdot x},\quad {{k}}\in \Gamma(\lambda),\quad x\in\mathbb{Z}^{d},
\label{discinc2}
\end{equation}
and $\psi_{sc}^{\rvs{-}}(x, k)$ is the outgoing solution for the non-homogenous equation
\begin{equation}
(\Delta+v)\psi_{sc}-\lambda\psi_{sc}=-v\psi_{0},
\label{discschro3}
\end{equation}
obtained using the limiting absorption principle, i.e. via formulas \eqref{psipm} and \eqref{asym1} with sign ${\rvs{-}}$, 
where $f=-v\psi_{0}$ and $\Delta$ is replaced by $\Delta+v$; see \cite{Shaban}
\rvs{and Remark \ref{remarkcontn}}.

We recall that, under assumption \eqref{assumeE1},
$\psi^{{\rvs{-}}}$ has the asymptotic expansion
\begin{eqnarray}
\psi^{{\rvs{-}}}(x, k) &=e^{i{{k}}\cdot{{x}}}+\dfrac{e^{{\rvs{-}}i\mu\left(\omega, {{\lambda}}\right) \left| {{x}}\right| }
}{\left| {{x}}\right| ^{\frac{d-1}2}}A\left({{k}},\omega\right)
+O(\dfrac1{\left| {{x}}\right| ^{\frac{d+1}2}})\nonumber \\
&\qquad\quad\textrm{as
}\left| {{x}}\right| \longrightarrow\infty,
\quad \omega=\frac{x}{|x|}, \quad x\in\mathbb{Z}^d,
\label{asympsiplus2}
\end{eqnarray}
where $\mu$ is as in \eqref{asym}, \eqref{defmu}
and the coefficient $A\left({k},\omega\right) $
is smooth and the remainder can be estimated uniformly in
$\omega$; see \cite{Eskina2,Shaban}.
The coefficient $A$ is the scattering amplitude for \eqref{discschro2}.

\begin{remark}
 In our previous article \cite{NovBls}, Eq. \eqref{discschro2} is misprinted with $-\Delta$ in place of $\Delta$;
see Eq. (1) and Eq. (10) in \cite{NovBls}.
However, there exists some disagreement regarding the standard conventions to 
present
the discrete Schr\"{o}dinger operators on $\mathbb{Z}^d$; compare, for example, the conventions of \cite{Shaban} and \cite{Isozaki1}.
\end{remark}

\begin{remark}
In the present article, we define $\mu$ exactly as in \cite{Shaban}.
This $\mu$ differs in sign from that defined in our previous article if $2d-4<\lambda<2d$. See also Remark 1.2 in \cite{NovBls}.
\end{remark}
In this appendix, we consider the inverse scattering problem for Eq. \eqref{discschro2} consisting in recovering the potential $v$ from the scattering amplitude $A$, for simplicity, under assumption \eqref{assumeE1}. 

As we mentioned already in introduction, many important results are given in the literature on inverse scattering for Eq. \eqref{discschro2}.
Below in this appendix,
we
discuss relations between inverse scattering for Eq. \eqref{discschro2} in the Born approximation for small $v$ and the inverse source problem for Eq. \eqref{discschro}, for simplicity, under assumption \eqref{assumeE1}. 

In the Born approximation for small $v$, under assumption \eqref{assumeE1}, we have that 
\begin{equation}
A(k,\omega)=-\frac{\sqrt{2\pi}\dF{v}\left(\newgamma\left({\rvs{-}} \omega,\lambda\right)-k \right) e^{\rvs{-}i(\sigma +2)\frac{\pi}{4}}}{\sqrt{\left| K(\omega,\lambda))\right| }\left| \nabla\phi(\newgamma\left(\omega,\lambda\right))\right| },\quad 
k\in\Gamma(\lambda),\quad
\omega\in \mathbb{S}^{d-1},
\label{ampshabnew}
\end{equation}
where $A$ is the scattering amplitude in \eqref{asympsiplus2},
$\dF{v}$ is the Fourier transform of $v$, defined according to \eqref{uFT},
$K$, $\phi$, $\newgamma, \sigma$ are the same as in \eqref{ampshab}.

Formula \eqref{ampshabnew} follows from \eqref{psiplus2}-\eqref{discschro3}, where in \eqref{discschro3} the product $v\psi_{sc}$ is replaced by $0$ in the Born approximation, and from formula \eqref{ampshab} with $f=-v\psi_{0}$.

In connection with formula \eqref{ampshabnew},
we also consider $\dF{v}$ on
\begin{equation}
\mathcal{B}(\lambda)
:=\{p=\newgamma-k : \newgamma, k\in\Gamma(\lambda)\},
\label{ballB}
\end{equation}
where
$\Gamma(\lambda)$
is as in Remark
\ref{remark1p1fromfirst},
and $\mathcal{B}(\lambda)$ can be considered as a subset of $\mathbb{R}^d$ as well as a subset of ${T}^d$.

\begin{remark}
    \rvs{For comparisions with the continuous case, in a way similar to Remark \ref{remarkcontn},
    note that $\Gamma(\lambda)\approx\{k\in\mathbb{R}^d: k^2=2d-\lambda\}, \phi(k)\approx 2d-k^2$ in a small neighbourhood of $\Gamma(\lambda)$,
    $\kappa(-\omega,\lambda)\approx \sqrt{2d-\lambda}~\omega$ as $\lambda\to 2d.$
    In view of these approximate equalities, the structure of \eqref{ampshabnew} is similar to the structure of formula for the scattering amplitude in the Born approximation for the continuous case; see, for e.g., \cite{Novikovrev}.
}
\end{remark}

The multi-frequency inverse source problem for Eq. \eqref{discschro}
and monochromatic inverse scattering for Eq. \eqref{discschro2} in the Born approximation for $d\ge2$ have considerable similarities following from the similarities of formulas \eqref{ampshab} and \eqref{ampshabnew}.
In particular, Theorems \ref{propAp1}, \ref{propAp2}, \rvs{and \ref{thmAstab}} given below for the latter case are analogues of Theorems \ref{thm1}, \ref{thm2}, \rvs{and \ref{thmstab}} for the former. 

\begin{theorem}
Let $v$ satisfy \eqref{defv2}, $\lambda$ satisfy \eqref{assumeE1}, and $d\ge2$.
Then the scattering amplitude $A$ for Eq. \eqref{discschro2}, arising in \eqref{ampshabnew} 
and given at fixed $\lambda$ in an open non-empty neighourhood $\mathcal{N}$ of any fixed pair $k, \omega,$
where $k\in\Gamma(\lambda)$ and $\omega\in\mathbb{S}^{d-1}$,
uniquely determines $v.$
\label{propAp1}
\end{theorem}

Under our assumptions, Theorem \ref{propAp1} follows from the statements:\\
(i) The function $A$ 
is real-analytic on $\Gamma(\lambda)\times \mathbb{S}^{d-1}$.
Therefore, $A$ on $\mathcal{N}$
uniquely determines $A$ on $\Gamma(\lambda)\times\mathbb{S}^{d-1}$.\\
(ii)
$A$ on $\Gamma(\lambda)\times\mathbb{S}^{d-1}$ uniquely determines $\dF{v}$ on $\mathcal{B}(\lambda)$ in \eqref{ballB} via \eqref{ampshabnew}.\\
(iii) $\dF{v}$ is real-analytic on ${T}^d.$
Therefore, $\dF{v}$ on $\mathcal{B}(\lambda)$
uniquely determines $\dF{v}$ on ${T}^d,$
and $v$ on $\mathbb{Z}^d$.

\begin{remark}
 Theorem \ref{propAp1} remains valid under assumption \eqref{assumeE1} but without the assumption that $v$ is small. 
 \rvs{For example, it remains valid under the assumptions of \cite{Isozaki2a}.}
 In this case, this result follows from the statement (i) as in the proof of Theorem \ref{propAp1} above
 and from the uniqueness for inverse scattering from full $A$ on $\Gamma(\lambda)\times\mathbb{S}^{d-1}$ proved in \cite{Isozaki2a}.
\end{remark}

\begin{remark}
    Similar to Remarks \ref{rem2p7}, \ref{rem2p8}, 
    our Theorem \ref{propAp1} remains valid for $v$ exponentially decaying at infinity, and is not valid, in general, for $v$ decaying at infinity as $O(|x|^{-\infty}).$
\end{remark}

\begin{theorem}
Let $u$ be a compactly supported function on $\mathbb{Z}^d$,
$\dF{u}={\mathcal{F}}u$,
and $\lambda$ satisfy \eqref{assumeE1}.
Let 
$v={\mathcal{F}}^{-1}\dF{v}$,
where
$\dF{v}(p)=\dF{u}(p)\prod_{j=1}^{J}(\phi(p+k_j)-\lambda)$, $\phi$ is defined by \eqref{defphi},
$p\in {T}^d$, $k_j\in\Gamma(\lambda)$, and $J$ is a positive integer.
Then
$v$ is compactly supported on $\mathbb{Z}^d$ and the scattering amplitude $A=A\left(k, \omega\right)$ arising in \eqref{ampshabnew} vanishes identically on $\mathbb{S}^{d-1}$ for each $k=k_j$, $j=1,\dotsc, J$.
\label{propAp2}
\end{theorem}

\begin{remark}
 The function $v$ in Theorem \ref{propAp2} can be also written as
 \begin{equation}
 v=\left(\prod\nolimits_{j=1}^J(L_{k_j} - \lambda) \right)u,
 \end{equation}
 where $L_{k_j} u(x)=e^{-ik_j x}\Delta(e^{ik_j x} u(x)),~ x\in\mathbb{Z}^d.$
 \label{remAp7}
\end{remark}

The proof of Theorem \ref{propAp2} is similar to the proof of Theorem \ref{thm2}.

\begin{remark}
Theorem \ref{propAp2} and Remark \ref{remAp7} admit straightforward polychromatic versions involving several $\lambda$.
\end{remark}

\begin{remark}
 Theorems \ref{propAp1} and \ref{propAp2} have analogues for the continuous case.
 In this case, Eq. \eqref{discschro2} is considered for $x\in\mathbb{R}^d$, where $\Delta$ is \rvs{replaced by the} the standard Laplacian $\triangle$,  $\phi(k)$ \rvs{replaced by} $\Phi(k)=-k^2, ~k\in\mathbb{R}^d$,
\rvs{and $\lambda<0$}. The analogue of Theorem \ref{propAp1} for the continuous case is well-known. 
 In the continuous analogue of Theorem \ref{propAp2}, the function $u$ should be also sufficiently smooth on $\mathbb{R}^d$ in order to ensure a regularity of the potential $v$ on $\mathbb{R}^d,$ and $\Gamma(\lambda)$ should be \rvs{replaced by} $\rvs{\mathbb{S}^{d-1}_{\sqrt{|\lambda|}}=}\{k\in\mathbb{R}^d: -k^2=\lambda\}$;
 see \cite{Devaney} for the continuous monochromatic case.
\end{remark}

\rvs{
\begin{theorem}
    Suppose that $v_1$ and $v_2$ satisfy \eqref{defv2},
$\lambda$ and $\mathcal{N}$ are the same as in Theorem \ref{propAp1},
and~ $\overline{\mathcal{N}}\subset\cup_{\zeta\in\Lambda} \Gamma(\zeta)\times\mathbb{S}^{d-1}$.
Then the following estimate holds:
    \begin{equation}
        \lVert v_2-v_1 \rVert_{\ell_2(D\cap\mathbb{Z}^d)}\le C_{D,\mathcal{N}} \lVert A_2-A_1 \rVert_{L_2(\mathcal{N})},
    \label{eqthmAstab}
    \end{equation}
    where
    $A_1$, $A_2$ are the scattering amplitudes arising in \eqref{ampshabnew} for $v_1, v_2$, respectively, and $C_{D,\mathcal{N}}$ is a positive constant depending on $D$ and $\mathcal{N}$ only.
\label{thmAstab}
\end{theorem}
}
\rvs{The proof of Theorem \ref{thmAstab} is similar to the proof of Theorem \ref{thmstab}.}

\rvs{
\begin{remark}
An open question consists in an extension of Theorem \ref{thmAstab} to the case when $v$ is not small.
In this connection,
one could combine techniques of \cite{Isozaki2a} and \cite{Bourgeois}.
     \label{conjAp12}
\end{remark}
}

\rvs{In connection with the case of phaseless inverse scattering},
Theorems \ref{propAp4} and \ref{propAp5}, given below 
for Eq. \eqref{discschro2} in Born approximation, are similar to Theorems \ref{thm5} and \ref{thm6} for phaseless inverse source problem for Eq. \eqref{discschro}.

\begin{theorem}
Let $v_0, v$ satisfy assumptions \eqref{eqsupp1}, \eqref{eqsupp} in place of $f_0, f$, where $v_0\not\equiv0$,
$\textrm{dist } (D, \, D_0) > \textrm{diam } \, D.$
Let $A_1$ be the scattering amplitude for Eq. \eqref{discschro2} arising in the Born approximation in \eqref{ampshabnew}, under assumption \eqref{assumeE1}, with $v$ replaced by $v+v_0$.
Then the intensity $|A_1|^2=|A_1\left(k, \omega\right)|^2$ given in an open nonempty neighbourhood $\mathcal{N}$ of any fixed pair $k, \omega,$ where $k\in\Gamma(\lambda)$ and $\omega\in\mathbb{S}^{d-1}$, uniquely determines $v$, assuming that $v_0$ is known a-priori.
\label{propAp4}
\end{theorem}

\begin{theorem}
Let $v_0, v$ satisfy assumptions \eqref{eqsupp1}, \eqref{eqsupp} in place of $f_0, f$, where $v_0\not\equiv0$,
$\textrm{dist } (D, \, D_0) > 0.$
Let $A$ and $A_1$ be the scattering amplitudes for $v$ and $v+v_0$, respectively,
arising 
in \eqref{ampshabnew}, 
under assumption \eqref{assumeE1}.
Then the intensities $|A|^2=|A\left(k, \omega\right)|^2$ and $|A_1|^2=|A_1\left(k, \omega\right)|^2$ given in an open nonempty neighbourhood $\mathcal{N}$ of any fixed pair $k, \omega$, where $k\in\Gamma(\lambda)$ and $\omega\in\mathbb{S}^{d-1}$, uniquely determines $v$, assuming that $v_0$ is known a-priori.
\label{propAp5}
\end{theorem}

The proofs of Theorems \ref{propAp4} and \ref{propAp5} are similar to the proofs of Theorems \ref{thm5} and \ref{thm6}.

\begin{remark}
    \rvs{
    In the present work, we consider reconstruction from far-field measurements.
    Reconstruction from near-field measurements can be reduced in many cases to
    reconstruction from far-field measurements.
    For discrete problems, results in this direction can be found in \cite{Isozaki2a,NovBls}.
    In particular,
    \cite{NovBls} gives formulas for finding phased far-field amplitude $A$ 
    from phaseless 
    measurements of $\psi^-$ for the discrete Schr\"{o}dinger equation \eqref{discschro2}.
    }
\end{remark}

\begin{remark}
\rvs{  Open questions include extensions of the results of 
the present work to the case of Helmholtz and Schr\"{o}dinger equations on more complicated lattices, for example, as in \cite{Ando,Bls6}.
}
\end{remark}

\begin{remark}
    \rvs{Open questions also include extensions of results of the present work to the case of inverse scattering for Eq. \eqref{discschro2}, and more general discrete Schr\"{o}dinger equations, in the framework of distorted Born approximation. In this connection for the continuous case, we refer, for example, to results of \cite{Chew,Novikov15,Hohage22}.}
\end{remark}

\end{document}